\documentclass[11pt]{article}
\usepackage{latexsym,amssymb,amsthm,amsmath,multirow}
\usepackage{amssymb}

\def\a{\alpha}\def\b{\beta }
\def\<{\left < }
\def\>{\right >}
\def\({\left ( }
\def\){\right )}
\def\e{\eqref}
\def\E4{\mathbb E^4 }

\def\mtc{\tilde M^n(4c)}

\def\ii{{\rm i}}
\def\tn{\tilde\nabla}

\def\de{\delta}

\def\n2{\left[{n\over2}\right]}
\font\b=cmr6
\font\b=cmr8

\begin{document}

\vbox{\hbox{ Proceedings of the Conference RIGA 2011}
\hbox{ {\it Riemannian Geometry and Applications}}
\hbox{ Bucharest, Romania}}
\vskip 1.5truecm

\setlength{\textheight}{19cm}
\setlength{\textwidth}{12.5cm}

\centerline{\large{\bf $\delta$-INVARIANTS FOR LAGRANGIAN SUBMANIFOLDS}}
\vskip.04in

\centerline{\large{\bf OF COMPLEX SPACE FORMS}}\medskip
\centerline{\bf Bang-Yen Chen and Franki Dillen}

\bigskip

\begin{abstract}  The famous Nash embedding theorem published in 1956 was aiming for the opportunity to use extrinsic help in the study of (intrinsic) Riemannian geometry, if Riemannian manifolds could be regarded as Riemannian submanifolds. However, this hope had not been materialized yet according to \cite{G}. The main reason for this was the lack of control of the extrinsic properties of the submanifolds by the known intrinsic invariants. In order to overcome such difficulties as well as to provide answers to an open question on minimal immersions, the first author introduced in the early 1990's new types of Riemannian invariants, his so-called $\delta$-curvatures, different in nature from the ``classical'' Ricci and scalar curvatures.

One purpose of this article is to present some old and recent results concerning $\delta$-invariants for Lagrangian submanifolds of complex space forms. Another purpose is to point out that the proof of Theorem 4.1 of \cite{CD11} is not correct and the Theorem has to be reformulated.  More precisely, Theorem 4.1 of \cite{CD11} shall be replaced by Theorems 8.1 and 8.3 of this article. Since the new formulation needs a new proof, we also provide the proofs of Theorems 8.1 and 8.3 in this article.

\vskip.04in


\end{abstract}

\section{Introduction}

Let $\tilde M^n$ be a complex $n$-dimensional  K\"ahler manifold endowed with the complex structure $J$ and the metric $g$.  The K\"ahler 2-form $\omega$ is defined by $\omega(\cdot\,,\cdot)=g(J\cdot,\cdot)$.  An isometric immersion  $\psi:M^n\to \tilde M^n(4c)$ of a Riemannian $n$-manifold $M^n$ into  $\tilde M^n$ is called {\it Lagrangian\/} if $\psi^*\omega=0$. Lagrangian submanifolds appear naturally in the context of classical mechanics and mathematical physics. For instance, the systems of partial differential equations of Hamilton-Jacobi type lead to the study of Lagrangian submanifolds and foliations in the cotangent bundle.

In differential geometry of submanifolds, theorems which relate intrinsic and extrinsic curvatures always play an important role. Related with the famous Nash embedding theorem \cite{nash},  the first author introduced in the early 1990's a new type of Riemannian invariants, denoted by $\delta(n_1,\ldots,n_k)$. He then established sharp general inequalities relating $\delta(n_1,\ldots,n_k)$ and the squared mean curvature $H^2$ for submanifolds in real space forms. Such invariants and inequalities have many nice applications to several areas in mathematics (see \cite{topics,book} for more details).
Similar inequalities also hold for Lagrangian submanifolds of complex space forms.

One purpose of this article is to present an incomplete survey of some old and recent results concerning $\delta$-invariants for Lagrangian submanifolds of complex space forms. Another purpose of this article is to point out that the proof of Theorem 4.1 of \cite{CD11} is not correct as stated and the theorem has to be reformulated.   More precisely, Theorem 4.1 of \cite{CD11} shall be replaced by Theorems 8.1 and  8.3 of this article. Since the new formulation needs a new proof, we also provide the proofs of Theorems 8.1 and 8.3 in this article.

\pagestyle{myheadings}

\markboth{B.-Y. Chen and F. Dillen}{ $\delta$-INVARIANTS}


\section{Preliminaries}

Let $\tilde M^n(4c)$ be a complete, simply-connected,  K\"ahler $n$-manifold  with constant holomorphic sectional curvature $4c$ and let $M^n$ be an $n$-dimensional  Lagrangian submanifold of  $\tilde M^n(4c)$. We denote the Levi-Civita connections of $M$ and  $\tilde M^n(4c)$ by $\nabla$ and $\tilde \nabla$, respectively.

The formulas of Gauss and Weingarten are given respectively by
\begin{align}\label{2.1} &\tn_X Y = \nabla_X Y + h(X,Y),\\ &\label{2.2} \tn_X \xi = -A_\xi X + D_X \xi,\end{align}
for tangent vector fields $X$ and $Y$ and normal vector fields  $\xi$, where $D$ is the normal connection. The second fundamental form $h$ is related to $A_\xi$ by
$$\<h(X,Y),\xi\> = \<A_\xi X,Y\>. $$
The mean curvature vector $\overrightarrow{H}$ of $M$ is defined by $$\overrightarrow{H}={1\over n}\,\hbox{trace}\,h.$$ The squared mean curvature $H^2$ is defined by
$H^2=\<\right.\! \overrightarrow{H},\overrightarrow{H}\!\left.\>.$

 For Lagrangian submanifolds, we have (cf. \cite{CO})
\begin{align}\label{2.3} &D_X JY = J \nabla_X Y,
\\\label{2.4} &A_{JX} Y = -J h(X,Y)=A_{JY}X.\end{align}
The above formulas immediately imply that $\<h(X,Y),JZ\>$ is totally symmetric.  If we denote the curvature tensors of $\nabla$ and $D$ by $R$ and $R^D$, respectively, then the equations of Gauss and Codazzi are given by
\begin{align} \label{2.5} & \<R(X,Y)Z,W\> =  \<A_{h(Y,Z)} X,W\> - \<A_{h(X,Z)}Y,W\>\\&
\notag \hskip1.1in +c(\<X,W\>\<Y,Z\>-\<X,Z\>\<Y,W\>),
\\ &\label{2.6} (\nabla h)(X,Y,Z) = (\nabla h)(Y,X,Z).\end{align}
where $X,Y,Z,W$ (respectively, $\eta$ and $\xi$) are vector fields
tangent (respectively, normal) to $M$; and $\nabla h$ is defined
by \begin{align}\label{2.7}(\nabla h)(X,Y,Z) = D_X h(Y,Z) - h(\nabla_X Y,Z) - h(Y,\nabla_X Z).\end{align}

For an orthonormal basis $\{e_1,\ldots,e_n\}$ of $T_pM$ at a point $p\in M$, we put $$h^A_{BC}=\<h(e_B,e_C),Je_A\>,\;\;A,B,C=1,\ldots,n.$$
It follows from \eqref{2.4} that \begin{align}
\label{2.8} h^A_{BC}=h^B_{AC}=h^C_{AB}.\end{align}


\section{Fundamental existence and uniqueness theorems}

For a Lagrangian submanifold in a K\"ahler manifold,
the cubic form $C$, defined by
$$C(X,Y,Z) = g(h(X,Y),JZ) = g(A_{JX}Y,Z),$$
is totally symmetric

The fundamental existence and uniqueness theorems for Lagrangian submanifolds are given by the following (cf. \cite{book,cdvv2}).

\vskip.1in
\noindent {\bf Theorem 3.1.} {\it Let $x^1,x^2 : M^n \to \tilde M^n(4c)$ be
two Lagrangian isometric immersions of a connected manifold $M^n$ into a complex space form $\tilde M^n (4c)$ of constant holomorphic sectional curvature $4c$.
If $$C^1(X,Y,Z) =C^2(X,Y,Z)$$
for all vector fields $X,Y,Z$ tangent to $M$, then there exists an
isometry $F$ of $\tilde M^n(4c)$ such that $x^1 = F(x^2)$.}

\vskip.1in
\noindent {\bf Theorem 3.2.} {\it  Let $(M^n,g)$ be an $n$-dimensional simply
connected Riemannian manifold. Let $\alpha$ be a symmetric bilinear
$TM^n$-valued form on $M^n$ satisfying
\begin{itemize}
\item $g(\alpha(X,Y),Z)$ is totally symmetric,
\item $(\nabla \alpha)(X,Y,Z)$ is totally symmetric,
\item $R(X,Y)Z = c X \wedge Y (Z) +
\alpha(\alpha(Y,Z),X)-\alpha(\alpha(X,Z),Y)$,
\end{itemize}
then there exists a Lagrangian immersion $x:M^n \to \tilde M^n(4c)$ such that
the second fundamental form $h$ satisfies $h(X,Y) = J \alpha(X,Y)$.}

Based on the fundamental existence theorem, we construct an important non totally geodesic Lagrangian immersion of a topological 3-sphere in complex projective space in the next chapter.


\section{An exotic Lagrangian immersion of $S^3$}

Consider the unit hypersphere $S^{2m+1}(1)\subset {\bf C}^{m+1}$ with standard Sasakian metric.
An immersion $f\colon\; M^n\to S^{2m+1}$ is called {\it C-totally real} or {\it horizontal} if for each
$p\in M$, ${\rm i} f(p)$ is normal to $M$. In particular, a $C$-totally real immersion $f\colon\; M^n\to S^{2n+1}$ is called {\it Legendrian}.

\medskip

Let $\pi : S^{2m+1}(1)\to CP^m(4)$ denote the Hopf fibration. If an immersion $f\colon\; M^n\to S^{2m+1}(1)$ is $C$-totally real (respectively, Legendrian), then $\pi(f)\colon\; M^n\to CP^m(4)$ is totally real (respectively, Lagrangian), and conversely (at least locally)  (cf. \cite{rec}).

Consider the unit sphere
$$
S^3 = \{(y_1,y_2,y_3,y_4) \in \mathbb R^4\,\vert\, y_1^2 + y_2^2+ y_3^2 +
y_4^2 = 1\}$$ in $\mathbb R^4$. Let $X_1$, $X_2$ and $X_3$ be the vector
fields defined by \begin{align*}
 &X_1(y_1,y_2,y_3,y_4) =
(y_2,-y_1,y_4,-y_3),\\
&X_2(y_1,y_2,y_3,y_4) = (y_3,-y_4,-y_1,y_2),\\
&X_3(y_1,y_2,y_3,y_4)= (y_4,y_3,-y_2,-y_1).\end{align*}
Let us define a metric $g$ on $S^3$ such that $X_1$, $X_2$ and $X_3$ are
orthogonal and
$$g(X_1,X_1)=
g(X_2,X_2) = 3,\,g(X_3,X_3)= 9.$$
We define a symmetric bilinear form $\alpha$ by
$$\begin{aligned}
&\alpha(X_1,X_1) = 2 X_1,\qquad&&\alpha(X_3,X_1)=0,\\
&\alpha(X_1,X_2) = -2 X_2,\qquad&&\alpha(X_3,X_2)=0,\\
&\alpha(X_2,X_2) = -2 X_1,\qquad&&\alpha(X_3,X_3)=0.
\end{aligned}
$$
Then $(S^3,g)$ and $\alpha$ satisfy all conditions of the
existence theorem for $c=1$. Hence there exists a Lagrangian immersion
$f: (S^3,g) \to CP^3(4).$
This Lagrangian immersion $f$ is minimal and $(S^3,g)$ has constant scalar curvature
$\frac13$, in particular $(S^3,g)$ is a {\it Berger sphere.}

An alternative description was given in \cite{BB} as follows. Define two complex structures on ${\bf C}^4$ by
\begin{equation}\begin{aligned}\notag
I(v_1,v_2,v_3,v_4)&=(iv_1,iv_2,iv_3,iv_4)\\
J(v_1,v_2,v_3,v_4)&=(-\bar v_4, \bar v_3,-\bar v_2,\bar v_1).
\end{aligned}\end{equation}
Clearly $I$ is the standard complex structure.
The corresponding Sasakian structures on $S^7(1)$ have characteristic vector
fields $\xi_1=-I(x)$ and $\xi_2=-J(x)$.
Since we consider two complex structures on ${\bf C}^4$, we can consider two different
Hopf fibrations $\pi_j:S^7(1)\to  C P^3(4)$. The vector field $\xi_j$ is vertical for $\pi_j$.

Now we consider the Calabi curve $\mathcal C_3$ of $C P^1$
into $C P^3(4)$ of constant Gauss curvature $4/3$, given by
$$
\mathcal C_3(z) = [1,\sqrt3 z ,\sqrt 3 z^2,z^3].
$$
Since $\mathcal C_3$ is holomorphic with respect to  $I$, there exists a circle
bundle $\pi : M^3 \to  CP^1$ over $CP^1$ and an isometric
minimal immersion $\mathcal I : M^3\to S^7(1)$ such that $\pi_1(\mathcal
I)=\mathcal C_3(\pi)$.  It is easy to check that $\mathcal I$ is horizontal
with respect to  $\pi_2$, such that the immersion $\mathcal J:M^3\to
CP^3(4)$, defined by $\mathcal J = \pi_2(\mathcal I)$, is a minimal Lagrangian isometric
immersion.  By straightforward computations, one
obtains that $\mathcal J$ has the required properties.
Thus this immersion is exactly the exotic immersion given above.

\section{Chen's $\delta$-invariants and fundamental inequalities.}

Let $M^n$ be an $n$-dimensional  Riemannian manifold. Denote by $K(\pi)$ the sectional curvature of $M$ associated with a plane section $\pi\subset T_pM^n$, $p\in M^n$. For any orthonormal basis
$e_1,\ldots,e_n$ of the tangent space $T_pM^n$, the scalar curvature $\tau$ at $p$ is non standardly
defined to be
\begin{align}\label{3.1}\tau(p)=\sum_{i<j} K(e_i\wedge e_j). \end{align}

More general, if $L$ is a subspace of $T_pM^n$  of dimension $r\geq 2$  and $\{e_1,\ldots,e_r\}$ an orthonormal basis of $L$, then the scalar curvature $\tau(L)$ of the $r$-plane section $L$ is defined by
\begin{align}\label{3.2}\tau(L)=\sum_{\alpha<\beta} K(e_\alpha\wedge e_\beta),\quad 1\leq \alpha,\beta\leq r.\end{align}

For given integers $n\geq 3$ and  $k\geq 1$, denote by $\mathcal S(n,k)$ the finite set  consisting of all $k$-tuples $(n_1,\ldots,n_k)$ of integers  satisfying  $$2\leq n_1,\cdots,
n_k<n\;\; {\rm and}\;\; n_1+\cdots+n_k\leq n.$$ Denote by ${\mathcal S}(n)$ the union $\cup_{k\geq 1}\mathcal S(n,k)$.

For each $(n_1,\ldots,n_k)\in \mathcal S(n)$ and each point $p\in M^n$, the first author introduced in \cite{c98,c00a} a Riemannian invariant $\de{(n_1,\ldots,n_k)}(p)$ defined  by
\begin{align}\label{3.3} \delta(n_1,\ldots,n_k)(p)=\tau(p)- \inf\{\tau(L_1)+\cdots+\tau(L_k)\},\end{align} where $L_1,\ldots,L_k$ run over all $k$ mutually orthogonal subspaces of $T_pM^n$ such that  $\dim L_j=n_j,\, j=1,\ldots,k$.

The invariants $\de(n_1,\ldots,n_k)$  and the  scalar curvature $\tau$ are very much different in nature (see \cite{topics} for a general survey on $\delta(n_1,\ldots,n_k)$).

Isometric Riemannian manifolds clearly have the same $\delta$-invariants. Therefore the $\delta$-invariants sometimes are called the DNA of the Riemannian manifold.

The first author proved in \cite{c98,c00a} the following optimal relationship between $\delta{(n_1,\ldots,n_k)}$ and the squared mean curvature $H^2$ for an arbitrary
submanifold in a real space form.

\vskip.1in
\noindent {\bf Theorem 5.1} {\it Let $M^n$ be an $n$-dimensional submanifold in a real space form $R^m(c)$ of constant  curvature $c$. Then, for each $k$-tuple  $(n_1,\ldots,n_k)\in\mathcal S(n)$, we have
\begin{equation}\begin{aligned}& \label{3.7} \delta{(n_1,\ldots,n_k)} \leq  {{n^2(n+k-1-\sum n_j)}\over{2(n+k-\sum n_j)}}H^2\\& \hskip1in +{1\over2} \Big({{n(n-1)}}-\sum_{j=1}^k {{n_j(n_j-1)}}\Big)c.\end{aligned}\end{equation}

The equality case of inequality \eqref{3.7}  holds at a point $p\in M$ if and only if, there exists an orthonormal basis  $\{e_1,\ldots,e_m\}$ at $p$, such that  the shape operators of $M$ in $R^m(\epsilon)$ at $p$  with respect to $\{e_1,\ldots,e_m\}$  take the form:
\begin{align} \font\b=cmr10 scaled \magstep2
\def\bigzerol{\smash{\hbox{ 0}}}
\def\bigzerou{\smash{\lower.0ex\hbox{\b 0}}} A_r=\left[ \begin{matrix} A^r_{1} & \hdots & 0
\\ \vdots  & \ddots& \vdots &\bigzerou \\ 0 &\hdots &A^r_k&
\\ \\&\bigzerou & &\mu_rI \end{matrix} \right],\quad  r=n+1,\ldots,m,
\label{3.8}\end{align}
where $I$ is an identity matrix and $A^r_j$ is a symmetric $n_j\times n_j$  submatrix satisfying
$$\hbox{\rm trace}\,(A^r_1)=\cdots=\hbox{\rm
trace}\,(A^r_k)=\mu_r.$$}
\vskip.1in

For $c=0$, the inequality \eqref{3.7} with $H=0$ can be considered as an obstruction for a Riemannian manifold the be immersible minimally in some Euclidean space, which gives a partial answer to a question of Chern.

The same result holds for Lagrangian submanifolds in a complex space form $\tilde M^n(4c)$ of constant holomorphic sectional curvature $4c$. More precisely, we have (cf. \cite{book}).

\vskip.1in
\noindent {\bf Theorem 5.2.} {\it  Let $M^n$ be an $n$-dimensional Lagrangian submanifold in a complex  space form $\tilde M^n(4c)$ of constant holomorphic sectional curvature $4c$. Then, for each $k$-tuple  $(n_1,\ldots,n_k)\in\mathcal S(n)$, we have
\begin{equation}\begin{aligned}\label{oldineq} \delta{(n_1,\ldots,n_k)} \leq &\, {{n^2(n+k-1-\sum n_j)}\over{2(n+k-\sum n_j)}}H^2\\& \hskip.1in +{1\over2} \Bigg({{n(n-1)}}-\sum_{j=1}^k {{n_j(n_j-1)}}\Bigg)c.\end{aligned}\end{equation}

The equality case of inequality \eqref{oldineq}  holds at a point $p\in M$ if and only if, there exists an  orthonormal basis  $\{e_1,\ldots,e_m\}$ at $p$, such that  the shape operators of $M$ in $\tilde M^n(4c)$ at $p$  with respect to $\{e_1,\ldots,e_m\}$  take the form of \eqref{3.8}.}
\vskip.1in

\vskip.1in

The following result was proved in \cite{c00b}.

\vskip.1in
\noindent {\bf Theorem 5.3.} {\it  Every Lagrangian submanifold of a complex space form $\tilde M^n(4c)$  that satisfies the equality case of inequality \eqref{oldineq} at a point $p$
 for some $k$-tuple $(n_1,\ldots,n_k)\in\mathcal S(n)$ is minimal at $p$}.

\vskip.1in
Theorem 3.3 extends  a result in \cite{cdvv1,cdvv2} on $\delta(2)$.


\section{The first Chen inequality}

A special case of Theorem 5.2 is for $k=1$ and $n_1=2$. In fact, the invariant $\delta(2)$ was introduced first and most results on $\delta$-invariants deal with $\delta(2)$. If we denote $\delta(2)$ by $\delta_M$, then Theorem 5.2 reduces to the following.
\vskip.1in

\noindent {\bf Theorem 6.1.} {\it
If $M^n$ is Lagrangian submanifold of $\mtc$, then
\begin{align}\label{chenfirst}\delta_M \leq  \frac{n^2(n-2)}{2(n-1)}H^2 + \tfrac12(n+1)(n-2)c.\end{align}

Equality holds at a point $p$ of $M$ if and only if the shape operators $A_r$ take the
following forms:
$$
A_{n+1} = \begin{pmatrix}
a&0&0&\dots&0\\
0&b&0&\dots&0\\
0&0&a+b&\dots&0\\
\vdots&\vdots&\vdots&\ddots&\vdots\\
0&0&0&\dots&a+b\end{pmatrix},\;
A_{r} = \begin{pmatrix}
h_{11}^r&h_{12}^r&0&\dots&0\\
h_{12}^r&-h_{11}^r&0&\dots&0\\
0&0&0&\dots&0\\
\vdots&\vdots&\vdots&\ddots&\vdots\\
0&0&0&\dots&0\end{pmatrix}.
$$}

For Lagrangian submanifolds satisfying the equality case of \e{chenfirst} we have the following results from \cite{cdvv2}.

\vskip.1in
\noindent {\bf Theorem 6.2.} {\it  Let $x:M^n \to \mtc$ ($n \ge 3$) be a Lagrangian
isometric immersion. If $M^n$ realizes equality in the first Chen
inequality \e{chenfirst} at a point $p$, then $M$ is minimal at $p$ and $h$ takes the form}
\begin{align*}
h(e_1,e_1) &= \lambda Je_1, \qquad h(e_1,e_2) = -\lambda Je_2\\
h(e_2,e_2) &= -\lambda Je_1,\qquad h(e_i,e_j) = 0.
\end{align*}

In \cite{cdvv1} the $\delta$-invariant $\delta_M$ is used to characterize the exotic immersion of $S^3$.

\vskip.1in
\noindent {\bf Theorem 6.3.} {\it  Let $x:M^n \to \mtc$ ($n\geq 3$), $c=-1,0,1$, be a Lagrangian
immersion with constant scalar curvature.  If $M^n$ realizes
equality in the first Chen inequality identically, then $M^n$ is totally
geodesic or $n=3$, $c=1$ and $x$ is congruent to the exotic immersion
$(S^3,g) \to CP^3$.}

In \cite{BSVW} and \cite{BSV} $3$-dimensional Lagrangian submanifolds of complex projective 3-space are classified.

\section{Oprea's improvement of the first Chen inequality}

In \cite{Op} Oprea improves the first Chen inequality as follows.

\vskip.1in
\noindent {\bf Theorem 7.1.} {\it Let $(M^n,g)$ be a Lagrangian submanifold of $\mtc$. Then we have}
$$\delta_M(p)\leq  \frac{n^2(2n-3)}{2(2n+3)}H^2+\frac12{(n-2)(n+1)}c.$$

\vskip.1in

It is shown in \cite{bo} that this inequality is sharp and that the constant on the right-hand side cannot be improved, by constructing an example for which equality is attained at one point.

The following theorem was also proved in \cite{bo}.

\vskip.1in
\noindent {\bf Theorem 7.1.} {\it Let $M^n$ be a Lagrangian submanifold of a $CP^n(4)$ attaining
equality in the improved first Chen inequality at every point. If $n\geq
4$, then $M$ is minimal.}

\vskip.1in

Non-minimal 3-dimensional Lagrangian submanifolds of $CP^3(4)$ that satisfy the equality case of the improved inequality were studied in \cite{BV}. In particular, it were proved in \cite{BV} that each such submanifold can be constructed from a certain minimal Lagrangian surface in $CP^2(4)$.

\section{Improved general inequalities}

In \cite{CD11} we improved the Chen inequality \e{oldineq}. Unfortunately
Theorem 4.1 of \cite{CD11} is not correct as stated and has to be reformulated. It shall be replaced by Theorems 8.1 and 8.3 of this section.
We will use the following convention concerning indices.
\begin{equation*}\begin{aligned}& \alpha_i,\beta_i,\gamma_i \in \Delta_i,\;\;  i,j\in \{1,\ldots, k\};  \\&r,s,t \in \Delta_{k+1}; \:\; u,v\in \{N+2,\ldots,n\};\;\;\\& A,B,C\in \{1,\ldots,n\},\end{aligned}\end{equation*}
where $\Delta_1=\{1, \ldots, n_1\}$, and for $1\leq i\leq k$
\begin{equation*}\Delta_i= \{n_1+\cdots+n_{i-1}+1,\ldots, n_1+\cdots+n_{i}\}.\end{equation*}
and
\begin{equation*}\Delta_{k+1}= \{N+1,\ldots, n \}.\end{equation*}

\vskip.1in
\noindent {\bf Theorem 8.1.}  {\it Let $M^n$ be a Lagrangian submanifold of a complex space form $\tilde M^n(4c)$. For a given $k$-tuple  $(n_1,\ldots,n_k)\in {\mathcal S}(n)$, we put $N=n_1+\cdots+n_k$ and $A= \sum_{i=1}^k (2+ n_{i})^{-1}$.  If $A\leq \frac{1}{3}$ and $N<n$, then we have
\begin{equation}\begin{aligned}\label{14.1} &\delta(n_1,\ldots,n_k) \leq  \text{$\dfrac{n^2 \big\{n- N+ 3k-1-6\, {\sum_{i=1}^k} (2+ n_{i})^{-1}\big\}}{2\big\{n-N +3k+ 2 -6\,{\sum_{i=1}^k} (2+ n_{i})^{-1}\big\}}$} H^2
\\& \hskip1in +\text{$\frac{1}{2}$}\Big\{n(n-1)-\text{$\sum$}_{i=1}^k n_i(n_i-1)\Big\}c.
\end{aligned}\end{equation}
The equality sign holds at a point $p\in M^n$ if and only if there is an orthonormal basis $\{e_1,\ldots,e_n\}$ at $p$ such that with respect to this basis the second fundamental form $h$ takes the following form
\begin{equation}\begin{aligned} \label{14.2} &h(e_{\alpha_i},e_{\beta_i})=\text{$\sum$}_{\gamma_i}\! h^{\gamma_i}_{\alpha_i \beta_i} Je_{\gamma_i}\!+\frac{3\delta_{\alpha_i\beta_i} }{2+n_i}\lambda Je_{N+1},
\\& h(e_{\alpha_i},e_{\alpha_j})=0, \;\;  \text{$\sum$}_{\alpha_i\in \Delta_i} h^{\gamma_i}_{\alpha_i\alpha_i}=0,
\\&h(e_{\alpha_i},e_{N+1})=\frac{3\lambda}{2+n_i} J e_{\alpha_i},\;\;   h(e_{\alpha_i},e_u)=0,
\\& h(e_{N+1},e_{N+1})=3\lambda Je_{N+1},\;
\; h(e_{N+1},e_u)=\lambda Je_u ,
\\&h(e_u,e_v)=\lambda \delta_{uv} Je_{N+1},\end{aligned}\end{equation} for distinct $ i,j=1,\ldots,k; \, u,v=N+2,\ldots,n ;$ and $ \lambda=\frac{1}{3}h^{N+1}_{N+1 N+1}$.}

\begin{proof} Let $\,(n_1,\ldots,n_k)\in \mathcal S(n)$ and let  $L_1,\ldots,L_k$ be mutually orthogonal subspaces of $T_pM$ with $\dim L_j=n_j$, $j=1,\ldots,k$.
We choose an orthonormal basis $\{e_1,\ldots,e_{n}\}$ at a point $p\in M$ which satisfies
$$e_{1},\ldots,e_{n_{1}}\in L_{1},\ldots, e_{n_1+\cdots+n_{k-1}+1},\ldots,e_{N}\in L_{k}.$$

Without loss of generality, we may assume that $c=0$.
Since
\begin{align} \label{14.3} &\tau=\sum_{A=1}^n\sum_{B<C} (h^A_{BB}h^A_{CC}-(h^A_{BC})^2),\\&\label{14.4}  \tau(L_i)=\sum_{A} \sum_{ \alpha_i<\beta_i} ( h^A_{\alpha_i \alpha_i}h^A_{\beta_i\beta_i}-(h^A_{\alpha_i\beta_i})^2),
\end{align}
we have
\begin{equation}\begin{aligned}\label{14.5}&\hskip-.1in
 \tau-\sum_{i=1}^k \tau(L_i)
 =\sum_{A}\sum_{r<s}(h^A_{rr}h^A_{ss}-(h^A_{rs})^2)+	\sum_{A,i}\sum_{\alpha_i,r} (h^A_{\alpha_i\alpha_i}h^A_{rr}- (h^A_{\alpha_i r})^2)
 \\& \hskip1.0in +\sum_A\sum_{i<j} \sum_{\alpha_i,\alpha_j}(h^A_{\alpha_i \alpha_i}h^A_{\alpha_j\alpha_j}-(h^A_{\alpha_i\alpha_j})^2)
\end{aligned}\end{equation}\begin{equation}\begin{aligned}\notag &\hskip.4in \leq \sum_{A}\Big\{\sum_{r<s} h^A_{rr}h^A_{ss}+	 \sum_{i}\sum_{\alpha_i,r} h^A_{\alpha_i\alpha_i}h^A_{rr}
 +\sum_{i<j} \sum_{\alpha_i,\alpha_j}h^A_{\alpha_i \alpha_i}h^A_{\alpha_j\alpha_j}\Big\}
\\& \hskip.9in
- \sum_{i}\sum_{\alpha_i,s}(h^{\alpha_i}_{ss})^2-\sum_{r\in\Delta_{k+1}}\sum_{B\ne r} (h^r_{BB})^2,
\end{aligned}\end{equation} with the equality sign holding if and only if
\begin{equation}\label{14.6}
h^{\alpha_i}_{\alpha_j \alpha_\ell} =h^{\alpha_j}_{\alpha_i\beta_i}=h^r_{\alpha_i \alpha_j}=h^{\alpha_i}_{st} =h^r_{st}=0 \end{equation}
for distinct $i,j,\ell\in\{1,\ldots,k\}$ and distinct $r,s,t\in \Delta_{k+1}$ and
\begin{equation}\label{14.6b}
h_{\alpha_i \beta_i}^r=0 \text{ for } \alpha_i \ne\beta_i.
\end{equation}

For a given $i\in \{1,\ldots,k\}$ and a given $\gamma_i\in \Delta_i$, we have
\begin{equation}\begin{aligned}\notag
 &0\leq \sum_{j=1}^k\sum_{r\in \Delta_{k+1}}\! \!\Big(\sum_{\alpha_j\in \Delta_j} h^{\gamma_i}_{\alpha_j \alpha_j}-3h^{\gamma_i}_{rr}\Big)^2+3\sum_{r<s}(h^{\gamma_i}_{rr}-h^{\gamma_i}_{ss})^2
\\&\hskip.1in + 3\sum_{\ell<j}\! \Big(\!\sum_{\alpha_\ell\in \Delta_\ell} h^{\gamma_i}_{\alpha_\ell \alpha_\ell}-\sum_{\alpha_j\in \Delta_j} h^{\gamma_i}_{\alpha_j \alpha_j}\Big)^2
\\& =(n\!-\! N\!+\!3k\!-\!3)\sum_{j}\Big(\sum_{\alpha_j} h^{\gamma_i}_{\alpha_j \alpha_j}\Big)^2\! -6\sum_{j}   \sum_{\alpha_j,r} h^{\gamma_i}_{\alpha_j \alpha_j} h^{\gamma_i}_{rr}-6 \sum_{r<s}h^{\gamma_i}_{rr}h^{\gamma_i}_{ss}
\\&\hskip.2in
 -6 \sum_{\ell<j} \sum_{\alpha_j}\sum_{\alpha_\ell} h^{\gamma_i}_{\alpha_\ell \alpha_\ell} h^{\gamma_i}_{\alpha_j \alpha_j} +3(n\!-\!N\!+\!3k\!-\!1)\sum_r (h^{\gamma_i}_{rr})^2
\\&=(n\!-\!N\!+\!3k\!-\!3)(h^{\gamma_i}_{11}+\cdots+h^{\gamma_i}_{nn})^2
- 2(n-N+3k)\times \\& \hskip.2in \Big\{\sum_{r<s}h^{\gamma_i}_{rr}h^{\gamma_i}_{ss} +\sum_{j=1}^k \sum_{\alpha_j,r} h^{\gamma_i}_{\alpha_j \alpha_j} h^{\gamma_i}_{rr} + \sum_{\ell<j}\sum_{\alpha_j,\alpha_\ell} h^{\gamma_i}_{\alpha_\ell \alpha_\ell}h^{\gamma_i}_{\alpha_j\alpha_j}-\sum_s (h^{\gamma_i}_{ss})^2\Big\}.
 \end{aligned}\end{equation}
Thus we find
\begin{equation}\begin{aligned} \label{14.7}&  \sum_{r<s}h^{\gamma_i}_{rr}h^{\gamma_i}_{ss} +\sum_{j=1}^k \sum_{\alpha_j,r} h^{\gamma_i}_{\alpha_j \alpha_j} h^{\gamma_i}_{rr} + \sum_{\ell<j}\sum_{\alpha_j,\alpha_\ell} h^{\gamma_i}_{\alpha_\ell \alpha_\ell}h^{\gamma_i}_{\alpha_j\alpha_j}-\sum_s (h^{\gamma_i}_{ss})^2
\\& \hskip.2in \leq \text{\small$\frac{n-N+3k-3}{2(n-N+3k)}$} (h^{\gamma_i}_{11}+\cdots+h^{\gamma_i}_{nn})^2.
\end{aligned}\end{equation}
with the equality holding if and only  if
\begin{equation}\begin{aligned} \label{14.8}&\text{$\sum$}_{\alpha_j\in \Delta_j}\! h^{\gamma_i}_{\alpha_j\alpha_j}=3 h^{\gamma_i}_{ss},\;\; j=1,\ldots,k,\;\;  s\in \Delta_{k+1}.\end{aligned}\end{equation}

Since $A\leq \frac{1}{3}$, we have
\begin{align}\label{14.9} \frac{n-N+3k-3}{n-N+3k}\leq \text{$\dfrac{n- N+ 3k-1-6A}{n-N+3k+ 2-6A}$}.\end{align}
Thus we get from \eqref{14.7} that
\begin{equation}\begin{aligned} \label{14.10}& \hskip-.2in \sum_{r<s}h^{\gamma_i}_{rr}h^{\gamma_i}_{ss} +\sum_{j=1}^k \sum_{\alpha_j,r} h^{\gamma_i}_{\alpha_j \alpha_j} h^{\gamma_i}_{rr} + \sum_{\ell<j} \sum_{\alpha_j,\alpha_\ell} h^{\gamma_i}_{\alpha_\ell \alpha_\ell}h^{\gamma_i}_{\alpha_j\alpha_j}-\sum_s (h^{\gamma_i}_{ss})^2
\\& \hskip.0in \leq \dfrac{n- N+ 3k-1-6A}{2\{n-N+3k+ 2-6A\}}\(\sum_{A=1}^nh^{\gamma_i}_{AA}\! \)^2.
\end{aligned}\end{equation}

If the equality holds in \e{14.10} with $k=1$, then $A<1/3$ and \e{14.9} is a strict inequality. Thus,  \e{14.6}, \eqref{14.6b} and \e{14.8}-\e{14.10} yield
$\sum_{\alpha_1\in \Delta_1}\! h^{\gamma_1}_{\alpha_1\alpha_1}=3 h^{\gamma_1}_{ss}=0,  s\in \Delta_{k+1}$.
If the equality holds in \e{14.10} with $k>1$,  it follows from \e{14.6}, \eqref{14.6b} and \e{14.8} that
\begin{equation}\begin{aligned} \label{14.11}&\text{$\sum$}_{\alpha_j\in \Delta_j}\! h^{\gamma_i}_{\alpha_j\alpha_j}=3 h^{\gamma_i}_{ss}=0,\;\; i=1,\ldots,k,\; \gamma_i\in \Delta_i \;  s\in \Delta_{k+1}.\end{aligned}\end{equation}
Thus, we have \e{14.11} for any $k\geq 1$. Conversely, it is easy to verify that \e{14.11} implies the equality case of \e{14.10}.

Let us put
$w=\text{\small$\frac{2}{3}$}\Big\{n-N+3k+2-\text{\small$\sum_{j=1}^k \frac{6}{2+n_j}\Big\}$}.$
Since
\begin{align}\notag \sum_{i=1}^k \frac{n_i}{2+n_i}=k-\sum_{i=1}^k \frac{2}{2+n_i},\;\; \; \sum_{j\ne i}\frac{n_j}{2+n_j}=k-\sum_{j}\frac{2}{2+n_j}-\frac{n_i}{2+n_i},\end{align}
we find  for each $t\in \{N+1,\ldots, n\}$  that
\begin{equation}\begin{aligned}\notag
 &0\leq \sum_i \sum_{r\ne t} \text{\small$ \frac{2+n_i}{3n_i}$}\Big(\sum_{\alpha_i} h^{t}_{\alpha_i \alpha_i}\! - \text{\small$ \frac{3n_i}{2+n_i}$}h^{t}_{rr}\Big)^2\\&\hskip.2in  +\sum_{i}\sum_{\alpha_i<\beta_i}\text{\small$\frac{w}{n_i}$} (h^{t}_{\alpha_i\alpha_i}\!-h^{t}_{\beta_i \beta_i})^2 \! +\!\sum_{\substack{r<s\\ r,s\ne t}} (h^t_{rr}\!-\! h^t_{ss})^2
\\&\hskip.1in + \sum_{i<j} \! \Bigg(\!  \text{\small$ \frac{\sqrt{(2+n_i)n_j}}{\sqrt{(2+n_j)n_i}}$} \sum_{\alpha_i } h^{t}_{\alpha_i \alpha_i}\! -\!  \text{\small$ \frac{\sqrt{(2+n_j)n_i}}{\sqrt{(2+n_i)n_j}}$}  \sum_{\alpha_j} h^{t}_{\alpha_j \alpha_j}\!\Bigg)^2
\\&\hskip.1in  +\text{\small$ \frac{1}{3}$} \sum_{r\ne t}(h^t_{tt}-3 h^t_{rr})^2
+\sum_i \text{\small$  \frac{n_i}{2+n_i}$} \Big( h^t_{tt}-\text{\small$  \frac{2+n_i}{n_i} $}\sum_{\alpha_i} h^t_{\alpha_i \alpha_i}\Big)^2
\\&= \sum_i \Big\{(n-N+2)\text{\small$ \frac{2+n_i}{3n_i}$}-\text{\small$\frac{w}{n_i}$} +\sum_{j\ne i}\text{\small$\frac{(2+n_i)n_j}{(2+n_j)n_i}$}\Big\}\Big(\sum_{\alpha_i}h^t_{\alpha_i\alpha_i}\Big)^2
\end{aligned}\end{equation}\begin{equation}\begin{aligned}\notag
& \hskip.1in  -2\sum_{r\ne t}\sum_{i}\sum_{\alpha_i}h^t_{\alpha_i\alpha_i}h^t_{rr}
+\Big\{n-N+1+\sum_i \text{\small$ \frac{3n_i}{2+n_i}$}\Big\}\sum_{r\ne t}(h^t_{rr})^2
\\&\hskip.1in +w \sum_i  \sum_{\alpha_i} (h^t_{\alpha_i \alpha_i})^2 -2\sum_{i<j}\sum_{\alpha_i,\alpha_j} h^t_{\alpha_i\alpha_i}h^t_{\alpha_j\alpha_j}-2\sum_{\substack{r<s\\ r,s\ne t}} h^t_{rr} h^t_{ss}
\\&\hskip.1in +\text{\small$ \Big\{\frac{n\!-\!N\!-\!1}{3}$}+\sum_i \text{\small$  \frac{n_i}{2+n_i}$}\Big\}(h^t_{tt})^2\! -2h^t_{tt}\sum_{r\ne t} h^t_{rr}\! -2h^t_{tt}\sum_i \sum_{\alpha_i} h^t_{\alpha_i \alpha_i}
\\&=\text{\small$\frac{1}{3}$} \Big\{n-N+3k-1-\sum_i \text{\small$  \frac{6}{2+n_i}$}\Big\}\sum_i\Big(\sum_{\alpha_i}h^t_{\alpha_i\alpha_i}\Big)^2  -2\sum_{r\ne t}\sum_{i,\alpha_i}h^t_{\alpha_i\alpha_i}h^t_{rr}
\\& \hskip.3in
+\Big\{n-N+3k+1-\sum_i \text{\small$ \frac{6}{2+n_i}$}\Big\}\sum_{r\ne t}(h^t_{rr})^2 -2h^t_{tt}\sum_i\!\sum_{\alpha_i} h^t_{\alpha_i \alpha_i}
\\
&\hskip.1in  -2\sum_{i<j} \sum_{\alpha_i \alpha_j} h^t_{\alpha_i\alpha_i}h^t_{\alpha_j\alpha_j}-2\sum_{\substack{r<s\\ r,s\ne t}} h^t_{rr} h^t_{ss} -2h^t_{tt}\sum_{r\ne t} h^t_{rr}+w \sum_i  \sum_{\alpha_i} (h^t_{\alpha_i \alpha_i})^2
 \\&\hskip.1in +\text{\small$\frac{1}{3}$} \Big\{n-N+3k-1-\sum_i \text{\small$  \frac{6}{2+n_i}$}\Big\}(h^t_{tt})^2
\\ &=w\left\{\text{\small$\dfrac{n\!-\! N\!+\! 3k\!-\!1-A}{2\{n\!-\! N\!+\!3k\!+\! 2\!-6A\}} $}\Big(\displaystyle{\text{\small$ \sum_{A=1}^n$}}h^{t}_{AA}\! \Big)^2-\sum_{r<s} h^t_{rr}h^t_{ss}\right.\\&\hskip.1in -
 \sum_{i<j} \sum_{\alpha_i,\alpha_j}\! h^t_{\alpha_i \alpha_i}h^t_{\alpha_j\alpha_j}  \!
 -\sum_{i}\sum_{\alpha_i,r} h^t_{\alpha_i\alpha_i}h^t_{rr}
+ \sum_{s\ne t} (h^t_{ss})^2 +\sum_i\sum_{\alpha_i} (h^t_{\alpha_i\alpha_i})^2\Bigg\}.
\end{aligned}\end{equation}
Hence, we obtain
\begin{equation}\begin{aligned} \label{14.12}&\sum_{r<s} h^t_{rr}h^t_{ss}+	\sum_{i}\sum_{\alpha_i,r} h^t_{\alpha_i\alpha_i}h^t_{rr}
 +\sum_{i<j} \sum_{\alpha_i,\alpha_j}h^t_{\alpha_i \alpha_i}h^t_{\alpha_j\alpha_j}
-\text{\small$ \sum_{B\ne t} $}(h^t_{BB})^2
\\& \hskip.4in \leq \dfrac{n- N+ 3k-1-6A}{2\{n-N+3k+ 2-6A\}} \(\sum_{A=1}^n h^{t}_{AA}\! \)^{\!2}.
\end{aligned}\end{equation}
with equality holding if and only  if  \begin{equation}\begin{aligned}\label{14.13} h^t_{tt}=(2+n_i)h^t_{\alpha_i\alpha_i}=3 h^t_{ss},\;\; i=1,\ldots,k,\;\; N+1\leq s\ne t\leq n.\end{aligned}\end{equation}
Thus, by combining \eqref{14.5}, \eqref{14.10} and \eqref{14.12}, we obtain inequality  \eqref{14.1}.

Equality in \eqref{14.1} implies that the inequalities \eqref{14.5},  \eqref{14.10} and  \eqref{14.12} become equalities. Thus, we have
\begin{align}\label{14.14} &\text{\small$\sum$}_{\alpha_j\in \Delta_j}\! h^{\gamma_i}_{\alpha_j\alpha_j}=3 h^{\gamma_i}_{ss}=0,\;\; j=1,\ldots,k,\;\;  s\in \Delta_{k+1},
\\&\label{14.15}h^t_{tt}=(2+n_i)h^t_{\alpha_i\alpha_i}=3 h^t_{ss},\;\; i=1,\ldots,k,\;\; N+1\leq s\ne t\leq n.
\\&\label{14.16} h^{\alpha_i}_{\alpha_j \alpha_\ell} =h^{\alpha_j}_{\alpha_i\beta_i}=h^r_{\alpha_i \alpha_j}=h^{\alpha_i}_{st} =h^r_{st}=0,\end{align}
for distinct $i,j,\ell\in\{1,\ldots,k\}$ and distinct $r,s,t\in \Delta_{k+1}$ and
\begin{equation}\label{14.16b}
h_{\alpha_i \beta_i}^r=0 \text{ for } \alpha_i \ne\beta_i.
\end{equation}
It follows from  \eqref{14.11} that the mean curvature vector lies in ${\rm Span}\{Je_{N+1},\ldots, $ $Je_n\}.$ Thus, we may choose $e_{N+1}$ in the direction of $JH$. So we conclude that conditions  \eqref{14.14}-\eqref{14.16b} are
 equivalent to  \eqref{14.2} due to the total symmetry of $h$.
\end{proof}

 The authors also proved  in \cite{CD11}  that the improved inequality given in Theorem 8.1 is best possible. More precisely, they proved the following.

\vskip.1in
\noindent {\bf Theorem 8.2.}  {\it
For each $k$-tuple $(n_1,\ldots,n_k)\in \mathcal S(n)$ satisfying $N<n$ and $A\leq \frac{1}{3}$, there exists a Lagrangian submanifold in ${\bf C}^n$ which satisfies the equality case at a point $p$ with $H(p)\ne 0$.}

\vskip.1in

Locally, every Lagrangian submanifold of ${\bf C}^n$ is given as graph:

\begin{equation*}
L(x_1,\ldots,x_n)=(x_1,\ldots,x_n, F_{x_1},\ldots, F_{x_n}),\;\; i=1,\ldots,n.
\end{equation*}
where $F=F(x_1,\dots,x_n)$ is any given function.
Take for $F$:
\begin{equation*}
F=\sum_{i=1}^k \frac{3\lambda }{2(2+n_i)}\sum_{\alpha_i\in \Delta_i} \!x_{\alpha_i}^2 x_{N+1}+ \frac{\lambda}{2}\sum_{r=N+1}^n\! x_{N+1} x_r^2.
\end{equation*}
It was shown in \cite{CD11} that the equality sign of \e{14.17} holds  at $0$ with $\overrightarrow{H}(0)\ne 0$ for the graph defined by this function $F$.

\vskip.1in
\noindent {\bf Theorem 8.3.}  {\it  Let $M^n$ be a Lagrangian submanifold of a complex space form $\tilde M^n(4c)$.   Then for any $k$-tuple  $(n_1,\ldots,n_k)\in {\mathcal S}(n)$ with $A> \frac{1}{3}$ and $N<n$ we have
\begin{equation}\begin{aligned}\label{14.17}  &\delta(n_1,\ldots,n_k) \leq  \dfrac{n^{2}(n- N + 3k-3)}{2(n-N+3k)} H^2
\\&\hskip.8in +\text{$\frac{1}{2}$}\Big\{n(n-1)-\text{$\sum$}_{i=1}^k n_i(n_i-1)\Big\}c.
\end{aligned}\end{equation}
The equality sign holds at a point $p\in M^n$ if and only if there exists an orthonormal basis $\{e_1,\ldots,e_n\}$ at $p$ such that
\begin{equation}\begin{aligned} \label{14.18} &h(e_{\alpha_i},e_{\beta_i})=\text{$\sum$}_{\gamma_i\in \Delta_{i}}\! h^{\gamma_i}_{\alpha_i \beta_i} Je_{\gamma_i},\;\;
\text{$\sum$}_{\alpha_i\in \Delta_{i}} h^{\gamma_i}_{\alpha_i\alpha_i}=0,
\;\;   \\& h(e_A,e_B)=0 \text{ otherwise},\end{aligned}\end{equation}
for $\alpha_{i},\beta_{i},\gamma_{i}\in \Delta_{i},\, i=1,\ldots,k$; and $A,B,C=1,\ldots,n$.}
\begin{proof}  Let $M^n$ be a Lagrangian submanifold of a complex space form.   Assume $A> \frac{1}{3}$.
Without loss of generality, we may assume that $c=0$. Then we have the following three inequalities as given in the proof of Theorem 8.1; namely,
\begin{equation}\begin{aligned}\label{14.19}&\hskip-.1in \tau-\sum_{i=1}^k \tau(L_i)
\leq \sum_{A}\! \Big\{\sum_{r<s} h^A_{rr}h^A_{ss}+\sum_{i}\sum_{\alpha_i,r} h^A_{\alpha_i\alpha_i}h^A_{rr}
 +\sum_{i<j} \sum_{\alpha_i,\alpha_j} \! h^A_{\alpha_i \alpha_i}h^A_{\alpha_j\alpha_j}\! \Big\}
\\& \hskip.9in
- \sum_{i}\sum_{\alpha_i,s}(h^{\alpha_i}_{ss})^2-\sum_{r=N+1}^n\sum_{B\ne r} \! (h^r_{BB})^2,
\end{aligned}\end{equation} with the equality sign holding if and only if \eqref{14.6} and \eqref{14.6b} hold;
\begin{equation}\begin{aligned} \label{14.21}&  \sum_{r<s}h^{\gamma_i}_{rr}h^{\gamma_i}_{ss} +\sum_{j=1}^k \sum_{\alpha_j,r} h^{\gamma_i}_{\alpha_j \alpha_j} h^{\gamma_i}_{rr} + \sum_{\ell<j} \sum_{\alpha_j \alpha_k}h^{\gamma_i}_{\alpha_\ell \alpha_\ell}h^{\gamma_i}_{\alpha_j\alpha_j}-\sum_s (h^{\gamma_i}_{ss})^2
\\& \hskip.2in \leq \frac{n-N+3k-3}{2(n-N+3k)} (h^{\gamma_i}_{11}+\cdots+h^{\gamma_i}_{nn})^2
\end{aligned}\end{equation}
with the equality holding if and only if
\begin{equation}\begin{aligned} \label{14.22}&\sum_{\alpha_j\in \Delta_j}\! h^{\gamma_i}_{\alpha_j\alpha_j}=3 h^{\gamma_i}_{ss},\;\; i,j =1,\ldots,k,\;\;  s\in \Delta_{k+1};\end{aligned}\end{equation}
and
\begin{equation}\begin{aligned} \label{14.23}&\sum_{r<s} h^t_{rr}h^t_{ss}+	\sum_{i}\sum_{\alpha_i,r} h^t_{\alpha_i\alpha_i}h^t_{rr}
 +\sum_{i<j} \sum_{\alpha_i,\alpha_j}h^t_{\alpha_i \alpha_i}h^t_{\alpha_j\alpha_j}
-\text{\small$ \sum_{B\ne t} $}(h^t_{BB})^2
\\& \hskip.4in \leq \text{\small$ \dfrac{n- N+ 3k-1-6A}{2\{n-N+3k+ 2-6A\}} $}\(\text{\small$\sum_{A=1}^n$} h^{t}_{AA}\! \)^{\!2}
\end{aligned}\end{equation}
with equality holding if and only  if  \begin{equation}\begin{aligned}\label{14.24} h^t_{tt}=(2+n_i)h^t_{\alpha_i\alpha_i}=3 h^t_{ss},\;\; i=1,\ldots,k,\;\; N+1\leq s\ne t\leq n.\end{aligned}\end{equation}

Since $A>\frac{1}{3}$, we have $k>1$ and
\begin{align}\label{14.25}\text{$\dfrac{n- N+ 3k-1-6A}{n-N+3k+ 2-6A}$}<\dfrac{n- N + 3k-3}{n-N+3k}.\end{align}
By combining \e{14.19}, \e{14.21}, \e{14.23} and \e{14.25} we obtain inequality \e{14.17}.

Now, let us assume that the equality sign of \e{14.17} at a point $p\in M^{n}$. Then it follows from
\e{14.22}, $k>1$, and $h^{\alpha_{j}}_{\alpha_{i}\beta_{i}}=0, i\ne j,$ in \e{14.6} that
\begin{equation}\begin{aligned} \label{14.26}&\sum_{\alpha_j\in \Delta_j}\! h^{\gamma_i}_{\alpha_j\alpha_j}= h^{\gamma_i}_{ss}=0,\;\; i,j =1,\ldots,k,\;\;  s\in \Delta_{k+1}.\end{aligned}\end{equation}
It follows from \e{14.17}, \e{14.23} and \e{14.25} that we also have $\sum_{A=1}^{n} h^{r}_{AA}=0$ for any $r\in \Delta_{k+1}$. Therefore, we conclude that $M^{n}$ is minimal in $\widetilde M^{n}(4c)$. Hence, from \e{14.6}, \e{14.6b}, \e{14.22}, \e{14.24} and \e{14.26}, and the minimality of $M^{n}$, we obtain \e{14.18}.

The converse can be verified directly.
\end{proof}

\vskip.1in
\noindent {\bf Remark 8.1.} Theorem 15.4 of \cite[p. 311]{book} and its proof shall be replaced accordingly.

\vskip.1in
\noindent {\bf Remark 8.2.}
Both inequalities improve the old inequality \eqref{oldineq}. For $A\leq 1/3$, the inequality cannot be improved. This follows from Theorem 8.2. For $A>1/3$ it is not known whether the inequality can be improved. So in Section 5 of \cite{CD11} the condition $A\leq 1/3$ should be assumed. On the other hand, the authors want to conjecture here that the condition $A\leq 1/3$ in Theorem 8.2. is not necessary. So we believe that Theorem 4.1 of \cite{CD11} is correct after all, but unfortunately we cannot prove \e{14.1} without assuming $A\leq 1/3$ yet.

\vskip.1in
\noindent {\bf Remark 8.3.}
Since $k=1$ automatically implies that $A<1/3$, the results in Section 6 and 7 of \cite{CD11} certainly remain valid.

\vskip.1in
\noindent {\bf Remark 8.4.}
In \cite{CD11} it is explicitly assumed that $N=n_1+\cdots+n_k<n$. Also in this section this is assumed. However, the $\delta$-invariant also makes sense if $N=n$ and the old inequality is also valid if $N=n$  (see \cite{c00a,book}), but it is not known whether it can be improved.


\section{Special case: $k=1$}

If  $k=1$, then $A=1/(n_1+2)\leq 1/4 < 1/3$. Thus we have the following result from Theorem 8.1.

\vskip.1in
\noindent {\bf Theorem 9.1.}  {\it   Let $M^n$ be an $n$-dimensional Lagrangian submanifold of a complex space form $\tilde M^n(4c)$. Then for any integer $n_1\in [2,n-1]$ we have
\begin{equation}\begin{aligned} &\delta(n_1)\leq \frac{n^2\{n_1(n-n_1)+2n-2\}}{2\{n_1(n-n_1)+2n+3n_1+4\}}H^2\\&\hskip.61in +\frac{1}{2}\{n(n-1)-n_1(n_1-1)\}c.\end{aligned}\end{equation}
Moreover, if  $M^n$ satisfies the equality case identically for some $n_1\leq n-2$, then $M^n$ is a minimal submanifold of $\tilde M^n(4c)$.}

\vskip.1in
When $n_1=n-1$, there exist non-minimal examples satisfying equality at every point (see \cite{CD11}).
All non-minimal Lagrangian submanifolds satisfying equality for $\delta(n-1)$ were classified recently by Chen, Dillen and Vrancken in \cite{cdv} as follows.

\vskip.1in
\noindent {\bf Theorem 9.2.}  {\it  Let $M^n$ be a non-minimal Lagrangian submanifold of the complex Euclidean $n$-space ${\mathbb{C}}^{n}$. Then
\begin{equation*} \delta(n-1)\leq \frac{1}{4}n(n-1)H^{2 }.\end{equation*}
Equality holds if and only if up to dilations and rigid motions, $M^n$ is defined by
\begin{equation*} L(\lambda, u_{2},\ldots,u_{n})=\frac{(n+1)e^{-\rm{i}\varphi}}{(n+1)\mu+\rm{i}\lambda}\phi( u_{2},\ldots,u_{n}), \end{equation*} where
\begin{equation*} \varphi(\lambda)=-\tfrac{n+1}{n}\csc^{-1}\!\((n+1)b \lambda^{\frac{n}{1-n}}\),\;\;
\mu(\lambda)=\sqrt{b^{2}\lambda^{\frac{2}{1-n}}-\tfrac{1}{(n+1)^{2}}\lambda^{2}}  \end{equation*}
for $b\in\mathbb{R}_0^+$, and $\phi(u_{2},\ldots,u_{n})$ is a minimal Legendrian submanifold of $S^{2n-1}(1)$.}

\vskip.1in
\noindent {\bf Theorem 9.3.}  {\it  Let $M^n$ be a non-minimal Lagrangian submanifold of the complex projective $n$-space $CP^{n}(4)$. Then
\begin{equation} \label{5.1} \delta(n-1)\leq \frac{n-1}{4}(nH^{2 }+4).\end{equation}
The equality sign holds identically if and only if  $\,M^n\,$ is congruent to the Lagrangian submanifold given by the composition $\pi\circ L$, where $\pi:S^{2n+1}(1)\to CP^{n}(4)$ is the Hopf fibration and
\begin{equation}\notag L(t, u_{2},\ldots,u_{n})=\( \frac{e^{-\ii\theta}\phi}{\sqrt{1+\mu^{2}+(n\!+\!1)^{-2}\lambda^{2}}},\frac{({\rm i} (n+1)^{{-1}}\lambda-\mu)e^{-n\ii \theta}}{\sqrt{1+\mu^{2}+(n\!+\!1)^{-2}\lambda^{2}}}\),  \end{equation}
where $\lambda(t), \mu(t)$ and $\theta(t)$  satisfy
\begin{equation}\begin{aligned} \notag  &\frac{d\theta}{dt}=-\frac{\lambda}{n+1},\;\; \lambda\ne 0,\\& \frac{d\lambda}{dt}=(n-1)\lambda\mu,\\& \frac{d\mu}{dt}=-1-\mu^{2}-\frac{n\lambda^{2}}{(n+1)^{2}}
\end{aligned} \end{equation}
and $\phi$ is a minimal Legendrian immersion in $S^{2n-1}(1)$.}

\vskip.1in
\noindent {\bf Theorem 9.4.}  {\it Let $M^n$ be a non-minimal Lagrangian submanifold of the complex hyperbolic $n$-space $CH^{n}(-4)$. Then
\begin{equation} \label{6.1} \delta(n-1)\leq \frac{n-1}{4}(nH^{2 }-4).\end{equation}
The equality sign holds identically if and only if $\,M^n\,$ is congruent to the Lagrangian submanifold given by the composition $\pi\circ L$, where $\pi:H^{2n+1}_{1}(-1)\to CH^{n}(-4)$ is the Hopf fibration and $L$ is one of the following immersions:

 {\rm (a)} $L(t, u_{2},\ldots,u_{n})=\( \frac{e^{-\ii\theta}\phi}{\sqrt{1-\mu^{2}-(n\!+\!1)^{-2}\lambda^{2}}},\frac{({\rm i}(n+1)^{{-1}}\lambda-\mu)e^{-n\ii \theta}}{\sqrt{1-\mu^{2}-(n\!+\!1)^{-2}\lambda^{2}}}\)$,
where $\lambda(t), \mu(t)$ and $\theta(t)$  satisfy
\begin{equation}\begin{aligned} \notag &\frac{d\theta}{dt}=-\frac{\lambda}{n+1},\;\; \lambda\ne 0,\\& \frac{d\lambda}{dt}=(n-1)\lambda\mu,\;\;
1>\mu^{2}+(n+1)^{-2}\lambda^{2},
\\& \frac{d\mu}{dt}=1-\mu^{2}-\frac{n\lambda^{2}}{(n+1)^{2}}
\end{aligned} \end{equation}
and $\phi$ is a minimal Legendrian immersion in $H^{2n-1}_{1}(-1)$;

{\rm (b)}   $L(t, u_{2},\ldots,u_{n})=\( \frac{({\rm i}(n+1)^{{-1}}\lambda-\mu)e^{-n\ii \theta}}{\sqrt{\mu^{2}+(n\!+\!1)^{-2}\lambda^{2}-1}},\frac{e^{-\ii\theta}\phi}{\sqrt{\mu^{2}+(n\!+\!1)^{-2}\lambda^{2}-1}}\)$,
where $\lambda(t), \mu(t)$ and $\theta(t)$  satisfy
\begin{equation}\begin{aligned}\notag  &\frac{d\theta}{dt}=-\frac{\lambda}{n+1},\;\; \lambda\ne 0,\\& \frac{d\lambda}{dt}=(n-1)\lambda\mu,\;\;
1<\mu^{2}+(n+1)^{-2}\lambda^{2},
\\& \frac{d\mu}{dt}=1-\mu^{2}-\frac{n\lambda^{2}}{(n+1)^{2}}
\end{aligned} \end{equation}
and $\phi$ is a minimal Legendrian immersion in $S^{2n-1}_{1}(1)$;

{\rm (c)}   the immersion given by
\begin{equation}\begin{aligned}\notag  &\hskip-.1in L(t, u_{2},\ldots,u_{n})=\frac{e^{\frac{2\ii }{n-1}\tan^{-1}(\tanh (\frac{1}{2}(n-1)t))}}{\cosh^{\frac{1}{n-1}}((n-1)t)} \times\\& \Bigg\{ \!\(w\!+\!\tfrac{\ii}{2}\!\<\phi,\phi\>\!+\!\ii,\phi,w\!+\!\tfrac{\ii}{2} \! \<\phi,\phi\>\)\\& \hskip-.2in +\left(\int_{0}^{t}\cosh^{\frac{2}{1-n}}((n-1)t)e^{2{\rm i}\tan^{-1}(\tanh (\frac{1}{2}(n-1)t))}dt\right)(1,0,\ldots,0,1)\Bigg\},
\end{aligned} \end{equation}
where $\phi(u_{2},\ldots, u_n)$ is a  minimal Lagrangian immersion in ${\bf C}^{n-1}$ and up to constants $w(u_2,\ldots,u_n)$ is the unique solution of the PDE system:
\begin{equation}\frac{\partial w}{\partial u_\a}=\left<{\rm i} \frac{\partial \phi}{\partial u_\a},\phi\right>,\;\; \a=2,\ldots,n.  \end{equation}}


\section{Pointwise symmetries}
Recall from Theorem 6.2 that the equality for $\delta_M$ in the minimal case holds if and only if
\begin{equation}
\begin{aligned}\label{S3}
h(e_1,e_1) &= \lambda Je_1, \qquad h(e_1,e_2) = -\lambda Je_2\\
h(e_2,e_2) &= -\lambda Je_1,\qquad h(e_i,e_j) = 0.
\end{aligned}
\end{equation}

Putting
$$f_1=\cos \! \(\frac{2\pi}{3}\)\! e_1 + \sin \!\(\frac{2\pi}{3} \)\! e_2, \, f_2=-\sin\! \(\frac{2\pi}{3}\)\! e_1 + \cos \! \(\frac{2\pi}{3}\)\! e_2, \, f_3=e_3,$$or
$$f_1=e_1, \, f_2=-e_2, \, f_3=-e_3,$$
doesn't change the form of \eqref{S3}.
Therefore $M$ has pointwise $S_3$-symmetry in the following sense.
For each $g\in S_3$, where $S_3$ is the group of order 6 generated by a rotation by angle $\frac{2\pi}{3}$ about the $z$-axis and a reflection in the $x$-axis one has, for all $u,v,w$ at the point $p$,
$$
C(gu,gv,gw)=C(u,v,w).
$$
This idea originates from Bryant's work.
One can consider similarly G-symmetry for any subgroup $G\subseteq \operatorname{SO}(n)$.

The following results are known.
\begin{enumerate}
\item For $n=3$ and $c=0$ the classification is done by R. L. Bryant \cite{B}.

\item For $n=3$ and $c\ne 0$, the classification can be deduced from work by L. Vrancken on affine spheres \cite{V}.

\item For $n=4$ and $c=0$, the classification is done by M. Ionel \cite{I}. But for $G=\operatorname{SO}(2)\rtimes S_3$, there are missing cases. In this case the second fundamental form takes the form \eqref{S3}.

\end{enumerate}

Recently, Dillen,  Scharlach,  Schoels and  Vrancken \cite{D} complete the classification for $G=\operatorname{SO}(2)\rtimes S_3$ in case $n=4$, and obtained a similar classification for any $c$.

\section{Remarks made on July 6, 2013.}

{\bf Remark 11.1.} Recently, it was proved in \cite{CDVV} that inequality (16) holds for arbitrary Lagrangian submanifolds in complex space forms (cf. Remark 8.2).
\vskip.1in

\noindent {\bf Remark 11.2.} A new general inequality with $N=n$ was proved in \cite{CDVV} for Lagrangian submanifolds in complex space forms (cf. Remark 8.4).

\hfill{\vbox{\hbox{B.-Y. Chen}
             \hbox{Department of Mathematics,}
             \hbox{Michigan State University,}
             \hbox{East Lansing, Michigan 48824, U.S.A.}
             \hbox{E-mail: {\tt bychen@math.msu.edu}}\hbox{}
             \hbox{F. Dillen}
             \hbox{Departement Wiskunde,}
             \hbox{Katholieke Universiteit  Leuven,}
             \hbox{Celestijnenlaan 200 B, Box 2400,}
             \hbox{BE-3001 Leuven (Belgium)}
             \hbox{E-mail: {\tt franki.dillen@wis.kuleuven.be}}}}

\end{document}